\documentclass[12pt,letterpaper]{article}
\usepackage{amsmath}
\usepackage{amsfonts}
\usepackage{algorithm}
\usepackage{verbatim}
\usepackage{geometry}  
\usepackage{graphicx}
\usepackage{caption}
\usepackage{mathtools}
\usepackage{multirow}
\usepackage{float}
\usepackage{amsthm}
\usepackage{amssymb}
\usepackage{epstopdf}
\usepackage{setspace}
\usepackage{fullpage}
\usepackage{enumerate}
\usepackage{authblk}
\usepackage{natbib}
\usepackage{hyperref}
\usepackage{bm}
\usepackage{color}

\newtheorem{definition}{\bf Definition}
\let\olddefinition\definition
\renewcommand{\definition}{\olddefinition\normalfont}

\newtheorem{theorem}{\bf Theorem}
\newtheorem{lemma}{\bf Fact}
\newtheorem{corollary}{\bf Corollary}

\DeclareMathOperator*{\argmax}{\arg\!\max}
\DeclareMathOperator*{\argmin}{\arg\!\min}

\newcommand{\bfx}{{\bf x}}
\newcommand{\umpbtg}{{\mbox{UMPBT}(\gamma)}}
\newcommand{\omgt}{{\Omega_\gamma(\theta)}}
\newcommand{\omgts}{{\Omega_\gamma(\theta^*)}}

\newcommand{\btt}{{\theta_t}}
\newcommand{\Hgamma}{{H_\gamma(\theta;\btt)}}
\newcommand{\con}{\,|\,}

\hypersetup{
	colorlinks,
	citecolor= blue
}

\graphicspath{{./Figures/}}

\begin{document}
	
\title{On the Existence of Uniformly Most Powerful Bayesian Tests With Application to Non-Central Chi-Squared Tests}

\author{Amir Nikooienejad and Valen E. Johnson\\
	\emph{Department of Statistics, Texas A\&M University}}

\date{}

\maketitle
		
		\begin{abstract}
			Uniformly most powerful Bayesian tests (UMPBT's) are an objective class of Bayesian hypothesis tests that can be considered the Bayesian counterpart of classical uniformly most powerful tests.  Because the rejection regions of UMPBT's can be matched to the rejection regions of classical uniformly most powerful tests (UMPTs), UMPBT's provide a mechanism for calibrating Bayesian evidence thresholds, Bayes factors, classical significance levels and p-values.  The purpose of this article is to expand the application of UMPBT's outside the class of exponential family models.  Specifically, we introduce sufficient conditions for the existence of UMPBT's and propose a unified approach for their derivation. An important application of our methodology is the extension of UMPBT's to testing whether the non-centrality parameter of a chi-squared distribution is zero.  The resulting tests have broad applicability, providing default alternative hypotheses to compute Bayes factors in, for example, Pearson's chi-squared test for goodness-of-fit, tests of independence in contingency tables, and likelihood ratio, score and Wald tests.  
		\end{abstract}
		
	
	\section{Introduction}\label{c5:intro}
	Bayesian hypothesis tests are based on computing the posterior probabilities of competing hypotheses given data. From Bayes theorem, the posterior probability of each hypothesis is proportional to the product of its prior probability and the marginal likelihood of the data given that the hypothesis is true.  In the case of two competing hypotheses, the posterior odds between hypotheses $H_0$ and $H_1$ can be written as 
	\begin{equation}\label{c5:postodds}
	\frac{\text{\bf{P}}(H_1\con\bfx)}{\text{\bf{P}}(H_0\con\bfx)} = \frac{m_1(\bfx)}{m_0(\bfx)} \times \frac{p(H_1)}{p(H_0)},
	\end{equation}
	where 
	$m_1(\bfx)/m_0(\bfx)$ is called the \emph{Bayes factor} in favor of the alternative hypothesis (denoted more simply as $\text{BF}_{10}(\bfx)$), $m_i(\bfx)$ denotes the marginal density of the data under hypothesis $i$, and $p(H_i)$ denotes the prior probability of hypothesis $H_i$. The logarithm of the Bayes factor is called the \emph{weight of evidence}. We assume throughout that the sampling density of the data $\bfx$ is defined with respect to a $\sigma-$finite measure and is described by the same parametric family of densities indexed by a parameter $\theta \in \mathbb{R}$ under all hypotheses, and refer to models and hypotheses interchangeably.  Letting $f(\bfx \con \theta)$ denote the sampling density of the data $\bfx$ given the value of a parameter $\theta \in \Theta$, and $\pi_i(\theta)$ the prior on $\theta$ given hypothesis $i$, the marginal density of the data under hypothesis  $i$ can be written as  
	\begin{equation*}
	m_i(\bfx) = \int_\Theta f(\bfx \con \theta)\pi_i(\theta) d\theta.
	\end{equation*}
	
	In the classical testing paradigm, a decision to reject the null hypothesis $H_0$ occurs when the value of a test statistic exceeds a specified threshold. 
 In a similar way, uniformly most powerful Bayesian tests (UMPBT's) of a null hypothesis are constructed by determining an alternative hypothesis that maximizes the probability that the Bayes factor of the test exceeds a pre-specified threshold, say $\gamma$, for all values of the data-generating parameter $\theta$.
	
	With this notation, a $\umpbtg$ was defined in \citet{umpbt} as follows:
	
	\begin{definition}\label{c5:firstdef}
		A uniformly most powerful Bayesian test for evidence threshold $\gamma>0$ in favor of the alternative hypothesis $H_1$ against a fixed null hypothesis $H_0$, denoted by UMPBT($\gamma$), is a Bayesian hypothesis test in which the Bayes factor for the test satisfies the following inequality for any $\theta_t \in \Theta$ and for all alternative hypotheses $H_2: \theta \sim \pi_2(\theta)$:
		
		\begin{equation}\label{c5:umpbtdef}
		\text{\bf{P}}_{\btt}\big[\text{BF}_{10}(\bfx) > \gamma\big] \geq \text{\bf{P}}_{\btt}\big[\text{BF}_{20}(\bfx) > \gamma\big].
		\end{equation}
	\end{definition}
	
	The alternative hypothesis $H_1$ in (\ref{c5:umpbtdef}) thus maximizes the probability that the Bayes factor is greater than a fixed evidence threshold, $\gamma$, among all possible prior densities that define alternatives hypotheses on the parameter space $\Theta$ and for all possible values of the data-generating parameter $\btt$ in $\Theta$. 
	
	For the case of testing simple null hypotheses $H_0:\theta=\theta_0$, and under the further assumption that tests are one-sided (i.e., $\Theta = \{\theta: \theta>\theta_0 \}$ or $\Theta = \{\theta: \theta<\theta_0 \}$), UMPBT's for one parameter exponential families were derived in \citet{umpbt}. These tests included tests of binomial proportions, tests of normal means with known variance, tests for normal variances when the mean is known, and tests that the non-centrality parameter of $\chi_1^2$ distribution is equal to zero \citep{pnas, umpbt}.  UMPBT's were extended in \citet{goddard} by restricting the class of alternative hypotheses over which the maximization in (\ref{c5:umpbtdef}) is performed.
	
	The UMPBT's derived in \citep{umpbt} were obtained by rewriting \newline
	$\text{\bf{P}}_{\btt}[\text{BF}_{10}(\bfx) > \gamma]$ in (\ref{c5:umpbtdef}) as
	\begin{equation}\label{c5:separable}
	\text{\bf{P}}_{\btt}\big[t(\bfx)> A(\gamma,\theta)\big],
	\end{equation}
	where $t(\bfx)$ was a sufficient statistic. By so doing, the probability in (\ref{c5:separable}) can be maximized with respect to $\theta$ by simply minimizing $A(\gamma,\theta)$, regardless of the distribution of $t(\bfx)$, thus producing a UMPBT($\gamma$) test. 
	
	The primary goal of this article is to provide a new approach to defining UMPBT's when $\text{\bf{P}}_{\btt}[\text{BF}_{10}(\bfx) > \gamma]$ cannot be written in the form of (\ref{c5:separable}). A primary application of the resulting method is to derive UMPBT's for tests of non-centrality parameters in $\chi^2$ distributions with arbitrary degrees of freedom, although the methodology presented in this article also requires the existence of scalar-valued sufficient statistics.  
	
	The extension of UMPBT's to $\chi^2$ tests is important because it facilitates the calibration of classical $p$-values and Bayes factors obtained from $\chi^2$ tests. This calibration can be accomplished by finding the evidence threshold $\gamma$ that produces the same rejection region as the classical test conducted at the given significance level $\alpha$.  This $\gamma$ implicitly defines the UMPBT alternative hypothesis, from which the Bayes factor of the test can be computed.  Thus, a correspondence between $\gamma$ and $\alpha$, and the $p$-value and the Bayes factor is obtained. Examples of this procedure are discussed in Section \ref{c5:discuss}

	Aside from applications involving hypothesis tests based on $\chi^2$ statistics, UMPBT's based on $\chi^2$ statistics have potential application in the realm of Bayesian model selection. For example, \citet{hu} propose the use of likelihood ratio test statistics to compute Bayes factors in model selection procedures.  For spurious regressors, the resulting $\chi^2$  test statistic has a central $\chi^2$ distribution; for important regressors it has a non-central $\chi^2$ distribution.  UMPBT's for the non-centrality parameter thus provide an objective prior for computing the marginal density of $\chi^2$-statistics in model selection procedures.  Importantly, the alternative model for the $\chi^2$ statistic implicit in this framework is a non-local prior density.   In contrast, other default Bayesian variable selection procedures are based on the use of local alternative priors on regression parameters \citep{intrinsic,ohagan}.  The potential value of using non-local priors in Bayesian variable selection is discussed in \citet{nlp}.

	The remainder of this article is organized as follows. Section \ref{c5:method} describes methodology to determine the existence of UMPBT's. In Section \ref{c5:commondist} we exploit this methodology to derive the $\umpbtg$ of a non-centrality parameter of a $\chi^2_\nu$ distribution with $\nu \geq 1$ degrees of freedom. Several examples are provided in Section \ref{c5:commondist}. 
	Concluding comments appear in Section \ref{c5:discuss}.
	
	\section{Method}\label{c5:method}
	\subsection{Preliminaries} 
	Let $y=h(\bfx)$ denote a sufficient statistic of the data, with $y \in \mathbb{R}$. For the remainder of this article, we assume that the null hypothesis being tested is a simple hypothesis having the form $H_0: \theta = \theta_0 \in \Theta$.
	
For every $\theta\in \Theta$, we denote the likelihood ratio in favor of $\theta_1$ as $g(y,\theta_1)$ (suppressing dependence on $\theta_0$).  For simple alternative hypotheses $H_1: \theta = \theta_1 \in \Theta$, $g(y,\theta_1)$ also represents the Bayes factor.    Our strategy for studying the properties of UMPBT's is to first restrict attention to Bayesian tests defined with simple alternative hypotheses, and to then extend properties of these tests to tests defined with composite alternative hypotheses.  If a UMPBT exists when the class of alternatives is restricted to simple alternatives, then the same test is a UMPBT when composite hypotheses are formed by averaging (according to a prior density) over simple alternative hypotheses.  This strategy is illustrated formally in the proof of Theorem 1.
	
For $\theta \in \Theta$, define 
	\begin{equation}\label{c5:rejreg}
	\omgt = \{y: g(y,\theta) > \gamma\}.
	\end{equation}
The interval $\Omega_\gamma(\theta_1)$ has a straightforward interpretation from a decision-theoretic perspective if we assume a $0-1$ loss function (i.e., the loss associated with incorrectly choosing the true hypothesis is 1, while the correct choice incurs no loss). In that case, $\Omega_\gamma(\theta_1)$ represents the ``rejection region'' of a simple null hypothesis against the simple alternative hypothesis $H_1: \theta=\theta_1$ when the prior odds in favor of $H_0$ are $\gamma$.
	
	Let $f(y;\btt)$ be the density function of $y$ for the data generating parameter, $\btt$, and $F(y;\btt)$ its corresponding distribution function defined with respect to a $\sigma$-finite measure $\mu$.  Let $S(f)\subset \mathbb{R}$ denote the support of $f$, which is assumed to be independent of $\theta$. Let $\bar{\mathbb{R}}$ represent affinely extended real numbers, $\mathbb{R}\cup\{-\infty,+\infty\}$. Define $a, b \in \bar{\mathbb{R}}$ as
	\begin{equation}\label{c5:ab}
	a=\inf S(f)\quad \quad b=\sup S(f).
	\end{equation}
	
	Next, define $H_\gamma(\theta_1;\theta_t) $ to be 
	\begin{equation}\label{c5:hfunc2}
	H_\gamma(\theta_1;\theta_t) = \text{\bf{P}}_{\btt}[g(y,\theta_1) > \gamma] = \int_{\Omega_\gamma(\theta_1)}f(y;\btt)\mu(dy),
	\end{equation}
	the probability that the Bayes factor exceeds $\gamma$ when the true state of nature is $\btt$ and the alternative is specified as $H_1:\theta = \theta_1$.
	
	If 
	\begin{equation}\label{c5:umpbt2}
	\theta^* = \argmax_{\theta \in \Theta} \Hgamma \quad \forall \btt \in \Theta,
	\end{equation}
	then it follows that $H_1:\theta=\theta^*$ is the alternative hypothesis of a $\umpbtg$.
	
The existence of UMPBT's, along with the fact that their alternatives do not place unit mass at the true parameter value, is somewhat counterintuitive. So, too, is the fact that they concentrate their mass on a (false) point alternative hypothesis. However, as noted in \citet{umpbt}, UMPBT's underestimate the ``true'' weight of evidence, (i.e., the logarithm of the Bayes factor) in favor of true alternative hypotheses in the sense that the expected value of the weight of evidence using the UMPBT alternative is less than it is under the true (data-generating) parameter value.  UMPBT's and the alternatives which define them thus provide a new class of default Bayesian hypothesis tests that are significantly less conservative than other default choices (e.g., \citet{jeffreys,intrinsic,ohagan,moreno98,bayarri08,bayarriforte}).
	
	\subsection{Existence and Derivation of UMPBT's}
	We now describe a sufficient condition for the existence of UMPBT's for one-sided hypothesis tests. The extension to two-sided tests requires further assumptions regarding the probability assigned to parameter values that are greater than or less than the null value.  If it is assumed that alternative hypothesis is symmetric around the null hypothesis, \citet{umpbt} showed that approximate two-sided UMPBT($\gamma$)'s can be obtained by specifying alternative hypotheses so that they concentrate their mass on the two corresponding one-sided UMPBT($2\gamma$) tests.
	
	\begin{theorem}\label{c5:coverthm}
		Consider a test of a simple null hypothesis $H_0:\theta=\theta_0$.
			Then $H_1: \theta = \theta^* \in \Theta$ defines an alternative hypothesis of a $\umpbtg$ if
			\begin{equation}\label{c5:covereq}
			\omgt \subseteq \omgts \quad \text{for all } \theta \in \Theta.
			\end{equation}
	\end{theorem}
	
	\emph{Proof:}
 For any simple alternative hypothesis $H_1: \theta\in \Theta$, the relation in (\ref{c5:covereq}) and the definition of the function $\Hgamma$ in (\ref{c5:hfunc2}) implies that
	\begin{equation}\label{c5:covh}
	\Hgamma = \int_{\omgt}F(dy;\btt) \leq \int_{\omgts}F(dy;\btt) = H_\gamma(\theta^*;\btt). 
	\end{equation}
	Knowing that $\theta^* \in \Theta$, the inequality above ensures that $\theta^* = \argmax_{\theta\in\Theta}\Hgamma$.  Because this inequality holds for all simple alternatives, it is straightforward to show that it also holds for composite alternatives \citep{umpbt}.  Let $\pi(\theta)$ be any prior density used to define the alternative hypothesis.  
		Define
		\begin{equation}
		s(y,\theta) = \begin{cases}
		1 \quad \mbox{if} \quad y \in \omgt \\
		0 \quad \mbox{otherwise.}
		\end{cases}
		\end{equation}
		Then (\ref{c5:covereq}) implies 
		\begin{equation} \int_{\Theta} s(y,\theta) \pi(d\theta) \leq 
		\int_{\Theta} s(y,\theta^*)  \pi(d\theta) = s(y,\theta^*) .
		\end{equation} 
		It follows that the Bayes factor based on $H_1: \theta \sim \pi(\theta)$, satisfies
		\begin{eqnarray}
		{\bf P}_{\theta_t}[BF_{10}(y)>\gamma] &=& \int_{\Theta}  \int_{\cal R} s(y,\theta) F(dy,\btt) \pi(d\theta)  \\
		&=&  \int_{\cal R} \int_{\Theta}  s(y,\theta) \pi(d\theta) F(dy,\btt)  \\
		&\leq &  \int_{\cal R}   s(y,\theta^*)  F(dy,\btt) = H_\gamma(\theta^*;\btt),
		\end{eqnarray}
		and the proof is complete. $\square$
	
	This is a useful existence theorem for UMPBT's. When $\omgt$ is an interval for all values of $\theta\in \Theta$, a more practical mechanism for establishing the existence of a UMPBT is provided in the following corollary.
	
	\begin{corollary}\label{c5:sufthm}
		Consider a Bayesian hypothesis test of a simple null hypothesis $H_0:\theta=\theta_0$.  
			If $\omgt$ is either of the form of $\big(c,d(\theta)\big)$ or $\big(d(\theta),c\big)$ for all $\theta \in \Theta$ and $c \in \overline{\mathbb{R}}$, and if 
			
			\begin{equation}\label{c5:thmeq}
			\theta^*=
			\argmin\limits_\theta vd(\theta), \; \mbox{ where }\; v=\begin{cases}
			-1 \; & \mbox{ if }\; \omgt = \big(c,d(\theta)\big)\\
			1 \; & \mbox{ if }\; \omgt = \big(d(\theta),c\big)
			\end{cases},
			\end{equation}
			then $H_1: \theta = \theta^*$ provides an alternative hypothesis corresponding to a UMPBT($\gamma$).
	\end{corollary}
	
	\emph{Proof:}
	It suffices to show that the condition (\ref{c5:covereq}) holds for the proposed $\theta^*$. Consider the case where $\omgt$ have the the form $\big(d(\theta),c\big)$. In this case, the upper bound of $\omgt$ is fixed for all $\theta \in \Theta$, but $\omgts$  has the smallest lower bound among all other $\omgt$. Hence,
	\begin{equation}\label{coreq}
	\omgt  \subseteq \omgts  \quad \mbox{for every } \theta\in\Theta .
	\end{equation}
	The proof follows from Theorem \ref{c5:coverthm}. The case for $v=-1$ follows along similar lines. $\square$
	
	Corollary \ref{c5:sufthm} offers a simple tool to check the existence of a UMPBT for continuous distributions, as well as offering a practical approach for finding it.
	
	If the Bayes factor is a monotone function of the sufficient statistic, then the following theorem provides a more direct route for finding a UMPBT.
	
	\begin{theorem}\label{c5:moncor}
		Let the likelihood ratio $g(y,\theta)$ be a continuous and differentiable function in $(a,b)\times\Theta$, the domains of $y$ and $\theta$. Define 
			\[ Q(\theta;y) = \frac{\partial g(y,\theta)}{\partial y},\] 
			and suppose for all $y$ and $\theta$ that $Q(\theta;y)$ is either strictly positive or strictly negative. Let $v$ denote the sign of $Q(\theta;y)$.
			For a fixed $\gamma>1$, let $\Lambda=\big\{(y,\theta): g(y,\theta)-\gamma = 0\big\}$ and suppose $\Lambda\neq\emptyset$. Let $r:\Theta\mapsto \mathbb{R}$ denote a function of $\theta$ such that $\big(r(\theta),\theta\big) \in \Lambda$. If 
			\begin{equation}\label{c5:moncorform}
			\theta^* = \argmin\limits_\theta vr(\theta),
			\end{equation} 
			then $H_1:\theta = \theta^*$ defines the alternative hypothesis for a UMPBT($\gamma$).
	\end{theorem}
	
	\emph{Proof:}
	Because $Q$ is either strictly positive or negative, the function $g$ is a one-to-one function of $y$.  Hence, for a given $\theta$, $g(y,\theta)-\gamma$ has a unique root. Due to the monotonicity of $g(y,\theta)$ in $y$, $\omgt$ is then either on the right side of the root, $\omgt=\big(r(\theta),b\big)$, or on its left, $\omgt=\big(a,r(\theta)\big)$, where $a$ and $b$ are defined in (\ref{c5:ab}). The proof follows from Corollary \ref{c5:sufthm}. $\square$ 
	
	Note that the form of $\omgt$ in the above theorem depends on $v$. More specifically, $\omgt$ is of the form $\big(r(\theta),b\big)$ when $v=1$ and it is of the form $\big(a,r(\theta)\big)$ when $v=-1$.  We note that for each value of $\theta$, the function $r(\theta)$ provides the value of $y=r(\theta)$ satisfying $g(y,\theta)=\gamma$; that is, the value of $y$ that results in a Bayes factor exactly equal to $\gamma$.  Because $Q$ is monotone, this value is unique.
		
		An example plot of $r(\theta)$ versus $\theta$ for different values of the evidence threshold, $\gamma$, is drawn in Supplementary Material for a non-central $\chi^2$ test with 10 degrees of freedom.
	
	Corollary \ref{c5:sufthm} and Theorem \ref{c5:moncor} provide sufficient conditions for the existence of a UMPBT. In Section \ref{c5:commondist} we apply these results to demonstrate both the existence of UMPBT's for one-parameter exponential family models and a UMPBT for the non-centrality parameter for $\chi^2$ distributions.
	
	Defining general conditions that are necessary for the existence of a UMPBT is difficult, but the next fact provides a simple method to demonstrate that a UMPBT does not exist.

		\begin{lemma}\label{noexlem}
			If the value of $\theta^*$ that maximizes $\Hgamma$ in (\ref{c5:hfunc2}) is not a constant function of $\theta_t$, then a $\umpbtg$ does not exist. 	
		\end{lemma}

	To see that this statement holds, suppose that $\theta^*_1$ maximizes $H_\gamma(\theta,\theta_{t_1})$ and $\theta^*_2$ maximizes $H_\gamma(\theta,\theta_{t_2})$, with $\theta^*_1\neq\theta^*_2$. It follows that there is no $\theta^*$ that maximizes $\Hgamma$ in (\ref{c5:umpbt2}) for every $\theta_t$ and thus a $\umpbtg$ does not exist. 
	
	An application of this fact to show that UMPBT's do not exist for one sample $t$-tests is provided in Section \ref{c5:onest}.
	
	\section{UMPBT's for Common Hypothesis Tests}\label{c5:commondist}
	
	\subsection{One-Parameter Exponential Family Distributions}
	We first show how the theory of the previous section can be used to derive UMPBT's for one-parameter exponential family distributions.  We also demonstrate that the method proposed in \citet{umpbt} is a special case of Theorem \ref{c5:moncor}. 
	
	Using the notation in \citet{umpbt}, let $\bfx = \{x_1,x_2,\cdots,x_n\}$ denote a random sample from a one-parameter exponential family model indexed by $\theta$, and suppose interest focuses on testing a null hypothesis $H_0: \theta=\theta_0$. 
	Our goal is to determine the UMPBT for a fixed evidence threshold, $\gamma$. We parametrize the density function for the model as
	\begin{equation}\label{c5:expfam}
	f(x\con\theta) = h(x)\exp\big\{\eta(\theta)T(x)-A(\theta)\big\},
	\end{equation}
	where $h(x), A(\theta)$ and $\eta(\theta)$ are known functions and $T(x)$ is the sufficient statistic of the data. Let $y=\sum_{i=1}^nT(x_i)$. For a simple alternative hypothesis $H_1:\theta = \theta_1$, it follows that the Bayes factor in favor of the alternative hypothesis can be expressed as
	\begin{equation}\label{c5:bfexpfam}
	g(y,\theta_1)= \exp\big\{n\big(A(\theta_0)-A(\theta_1)\big)\big\}\exp\Big\{y\big(\eta(\theta_1)-\eta(\theta_0)\big)\Big\}.
	\end{equation}
	Consequently, the first derivative of the Bayes factor with respect to $y$ in (\ref{c5:bfexpfam}) can be written as
	\begin{equation}\label{c5:derbfexpfam}
	\frac{\partial g(y,\theta_1)}{\partial y}=\big[\eta(\theta_1)-\eta(\theta_0)\big]\exp\Big\{n\big(A(\theta_1)-A(\theta_0)\big)+y\big(\eta(\theta_1)-\eta(\theta_0)\big)\Big\}.
	\end{equation}
	If the function $\eta(\theta_1)$ is monotonic on $\Theta$, the derivative above does not change sign and is strictly positive or negative. Therefore, for a fixed threshold $\gamma$, $g(y,\theta_1)-\gamma$ has a unique root, given by
	\begin{equation}\label{c5:rootexpfam}
	y= \frac{\log(\gamma)+n\big(A(\theta_1)-A(\theta_0)\big)}{\eta(\theta_1)-\eta(\theta_0)}.
	\end{equation}
	
	Following Theorem \ref{c5:moncor}, if 
	\begin{equation}\label{c5:expfameq}
	\theta^* = \argmin\limits_{\theta\in\Theta} v\frac{\log(\gamma)+n\big(A(\theta)-A(\theta_0)\big)}{\eta(\theta)-\eta(\theta_0)},
	\end{equation}
	where $v$ is equal to the sign of $\eta(\theta)-\eta(\theta_0)$, then $H_1:\theta=\theta^*$ defines an alternative hypothesis for a UMPBT($\gamma)$.
	
	In testing a one sided alternative against a point null hypothesis for one dimensional exponential family distributions, the $\umpbtg$ can always be found as described in (\ref{c5:expfameq}) if the natural parameter of the exponential family, $\eta(\theta)$, is monotone on the domain of the alternative hypothesis $\Theta$. This result  confirms the findings in \citet{umpbt}.
	
	\subsection{UMPBT's for the Non-centrality Parameter in $\chi^2$ Tests}
	We now apply the theory of Section 2 to derive UMPBT's for the non-centrality parameter of $\chi^2$ test statistics.  We then apply these tests to contingency tables and use them to study the relationship between $p$-values based on $\chi^2$ tests and Bayes factors obtained from the corresponding UMPBT.
	
	Let $y$ be an observation from a chi-squared distribution on $\nu$ degrees of freedom and non-centrality parameter $\theta$, denoted by $\chi^2_\nu(\theta)$. As shown in \citet{patnaik} and \citet{seber}, the probability density function of $y$ can be written as 
	\begin{equation}\label{c5:ncpdf}
	f(y\con\theta) = \frac{1}{2}\exp\left[-\frac{(y+\theta)}{2}\right]\big(\frac{y}{\theta}\big)^{\nu/4-1/2}I_{\nu/2-1}(\sqrt{\theta y}).
	\end{equation}
	Here, $I_\nu(\cdot)$ is the modified Bessel function of the first kind and for a real valued $\nu$ is defined as
	\begin{equation}\label{c5:besselfirst}
	I_\nu(y) = \sum\limits_{j=0}^\infty \frac{(y/2)^{2j+\nu}}{\Gamma(\nu+j+1)j!}.
	\end{equation} 
	
	In general, the range of the modified Bessel function of the first kind is $\mathbb{C}$, the set of all complex numbers. However, for real positive arguments and real-valued degrees of freedom, the range is $\mathbb{R}^+$. In the case of $\theta=0$, the probability distribution function in (\ref{c5:ncpdf}) reduces to
	\begin{equation}\label{c5:cpdf}
	f(y|\theta=0) =\left(\frac{1}{2}\right)^{\nu/2}e^{-y/2}\, \frac{y^{\nu/2-1}}{\Gamma(\nu/2)}.
	\end{equation}
	
	Suppose we are interested in testing $H_0: \theta=0$ against $H_1:\theta \sim \pi(\theta)$, where $\pi(\theta)$ is any probability density function defined on the non-negative real line not representing a point mass at 0. Using (\ref{c5:ncpdf}) and (\ref{c5:cpdf}), the Bayes factor in favor of a simple alternative hypothesis $H_1: \theta=\theta_1 \neq 0$ can be expressed as
	
		\begin{equation}\label{c5:bfncp}
		g(y,\theta_1)= \Gamma\left(\frac{\nu}{2}\right)\exp^{-\theta_1/2}2^{\nu/2-1}
		(\sqrt{\theta_1 y})^{1-\nu/2}I_{\nu/2-1}(\sqrt{\theta_1 y}).
		\end{equation}

	For this Bayes factor, both the data and the parameter of the test are arguments of the modified Bessel function. Thus ${\bf P}_{\theta_t}[BF_{10}(y)>\gamma]$ cannot be written in the form of (\ref{c5:separable}). The following theorem proves the existence of a $\umpbtg$ for this test using Corollary \ref{c5:sufthm}. 
	
		\begin{corollary}\label{c5:noncxthm}
			Suppose $y\sim \chi^2_\nu(\theta)$ and consider a test of the null hypothesis $H_0: \theta = 0$. 
			Given an evidence threshold $\gamma >0$, define $r(\theta)$ as in Theorem \ref{c5:moncor}. Then the alternative hypothesis that defines the $\umpbtg$ is given by $H_1: \theta = \theta^*$, where
			\begin{equation}\label{c5:nonceq}
			\theta^*=\argmin\limits_{\theta>0}r(\theta).
			\end{equation}
		\end{corollary}

	\emph{Proof:} The first derivative of the modified Bessel function of the first kind with $\nu$ degrees of freedom can be expressed as  $\frac{\partial I_\nu(z)}{\partial z}=\frac{\nu}{z}I_\nu(z)+I_{\nu+1}(z)$. The first derivative of $g(y,\theta)$ with respect to $y$ thus equals 
	\begin{equation}\label{c5:bfder}
	\frac{\partial g(y,\theta)}{\partial y} = \frac{\alpha}{2}\theta (\sqrt{\theta y})^{-\nu/2}I_{\nu/2}(\sqrt{\theta y}),
	\end{equation}
	where $\alpha = \Gamma(\frac{\nu}{2})\exp^{-\theta/2}2^{\nu/2-1}$ is a positive number. The domain for the alternative hypothesis is $\theta \geq 0$ and the support of the $\chi^2$ distribution is $\mathbb{R}^+$, which results in a real, positive modified Bessel function of the first kind. Therefore, the derivative in (\ref{c5:bfder}) is strictly positive. Moreover, $g(y,\theta)$ is continuous on $\Theta\times\mathbb{R}^+$ and its infimum is zero. Hence, for every $\gamma>0$, the set $\Lambda$ defined in Theorem \ref{c5:moncor} is not an empty set. The result  then follows from Theorem \ref{c5:moncor}. $\square$	
	
	\subsubsection{Tests of Independence in Contingency Tables}\label{c5:tct}
	Tests of independence between rows and columns of contingency tables are common in statistical practice. Performing such tests in the Bayesian paradigm requires computation of the Bayes factor, which depends on the prior densities assumed for the multinomial probability vector under both hypotheses. Different methods have been proposed to define these priors. \citet{albert} uses a prior distribution for the alternative hypothesis constructed about the ``independence surface'' representing the null hypothesis. \citet{good} used a mixed-Dirichlet prior and checked the robustness and sensitivity of this assumption with respect to hyperpriors and their hyperparameters. \citet{bfts} proposed a totally different approach by computing the Bayes factor based on a test statistic, in this case the standard $\chi^2$-statistic. Johnson's approach requires the specification of a prior distribution on the non-centrality parameter of the chi-squared distribution under the alternative hypothesis. He used a conjugate gamma prior density for the non-centrality parameter, and discussed various schemes for setting the hyperparameters of the prior density.
	
	Our method extends the concept of uniformly most powerful Bayesian tests to non-central chi-squared tests with different degrees of freedom. As a result, borrowing the methodology from \citet{bfts}, we use a $\chi^2$-statistic to compute the Bayes factor. The difference between our method and \citet{bfts} is that we use UMPBT methodology to fix the prior on the non-centrality parameter under the alternative hypothesis. 
	
	The contingency table shown in Table \ref{c5:ctable} represents the cross classification on cancer site and blood type for patients with stomach cancer \citep{white}. The total sample size is 707 and the goal is to test independence of cancer site and blood type.
	
	\begin{table}[!t]
		\centering	
		\caption{\citet{white} classification of cancer patients}
		\begin{tabular}{|lrrr|}
			\hline
			\multirow{2}{*}{Site} & \multicolumn{3}{c|}{Results for the} \\
			& \multicolumn{3}{c|}{following blood groups:}\\ \cline{2-4}
			& O & A & B or AB \\ \hline
			Pylorus and antrum & 104 & 140 & 52 \\
			Body and fundus & 116 & 117 & 52 \\
			Cardia & 28 & 39 & 11 \\
			Extensive & 28 & 12 & 8 \\ \hline		
		\end{tabular}
		\label{c5:ctable}
	\end{table}
	
	The $\chi^2$-statistic for this contingency table is $12.65$ on $6$ degrees of freedom. \citet{bfts} computed the Bayes factor against the independence model as $2.97$ when the hyperparameter of the prior gamma density was chosen so as to maximize the Bayes factor. 
	
	Using the method proposed by \citet{albert}, the maximum Bayes factor that can be obtained in favor of the alternative hypothesis is $3.02$. This value is obtained by maximizing the approximate Bayes factor with respect to the parameter that controls the dispersion of the alternative hypothesis around the independence model. Under the model proposed by \citet{good}, the Bayes factor is $3.06$. 
	
	Using the methodology proposed in this article, the Bayes factor based on the $\chi^2$-statistic is $3.52$. This value is obtained by assuming that $\omgts$  matches the rejection region of a 5\% classical test. The alternative hypothesis in this test is that the non-centrality parameter is equal to $7.31$. The evidence threshold corresponding to the 5\% test is $\gamma=3.46$.  
	
	In assessing the evidence against null hypotheses, this example illustrates that UMPBT's are not as conservative as other default Bayesian tests.  This is especially true when Bayesian tests specify local alternative hypotheses, or alternative hypotheses that place prior mass around the null value \citep{nlp}

	\subsubsection{Comparing Bayesian and Classical Tests of $\chi^2$ Non-Centrality Parameters }
	
\begin{figure}[t]
		\centering
		\includegraphics[width=75mm]{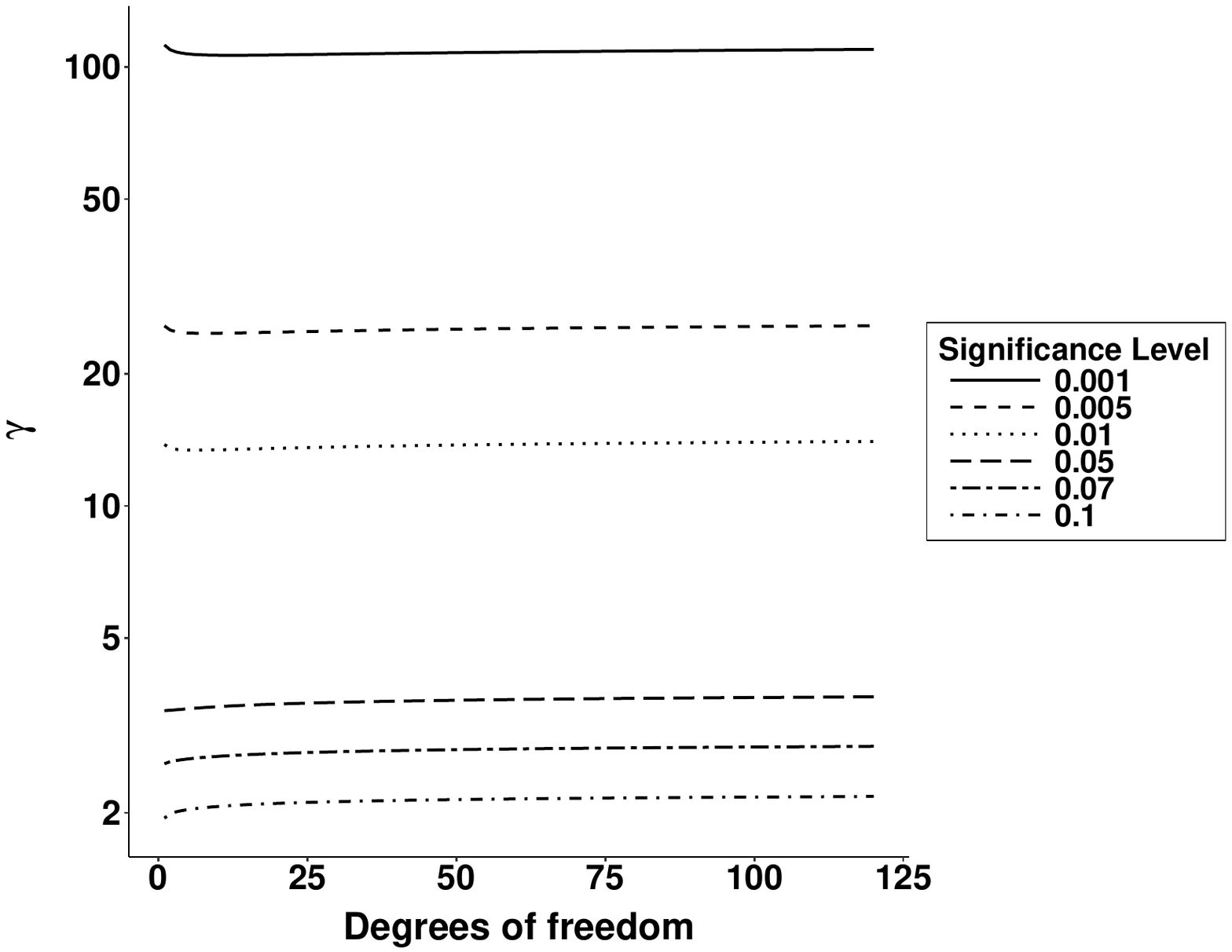}
		\caption{Evidence threshold vs. degrees of freedom in chi-squared tests for different significance thresholds.}
		\label{c5:thall}
	\end{figure}

	The connection between $p$-values and evidence thresholds based on $\chi^2$ tests for independence in contingency tables holds more broadly for hypothesis tests based on $\chi^2$ statistics. If $H_0: \theta_0=0$, and $H_1: \theta = \theta^*$ defines the alternative hypothesis of a UMPBT($\gamma$), where $\gamma$ is chosen so that $\omgts$ matches the rejection region of a classical  $\chi^2$ test of size $\alpha$, it is possible to compare evidence thresholds $\gamma$ to $p$-values in general tests of non-centrality parameters.  For instance, Figure \ref{c5:thall} illustrates the relation between the Bayesian evidence thresholds $\gamma$ and the p-values of classical tests (with rejection region matched to $\omgts$) versus the degrees of freedom of the $\chi^2$ tests. 
	
	Two important points are exposed in this figure. First, for a given p-value the evidence thresholds from corresponding Bayesian tests are almost constant with respect to the degrees of freedom. Second, for degrees of freedom less than 120, a $p$-value of 0.05 is equivalent to evidence thresholds that are always less than 3.67.  This value of the evidence threshold suggests positive, but not strong evidence, against the null hypothesis. 

	\begin{figure}[!t]
		\centering
		\includegraphics[width=75mm]{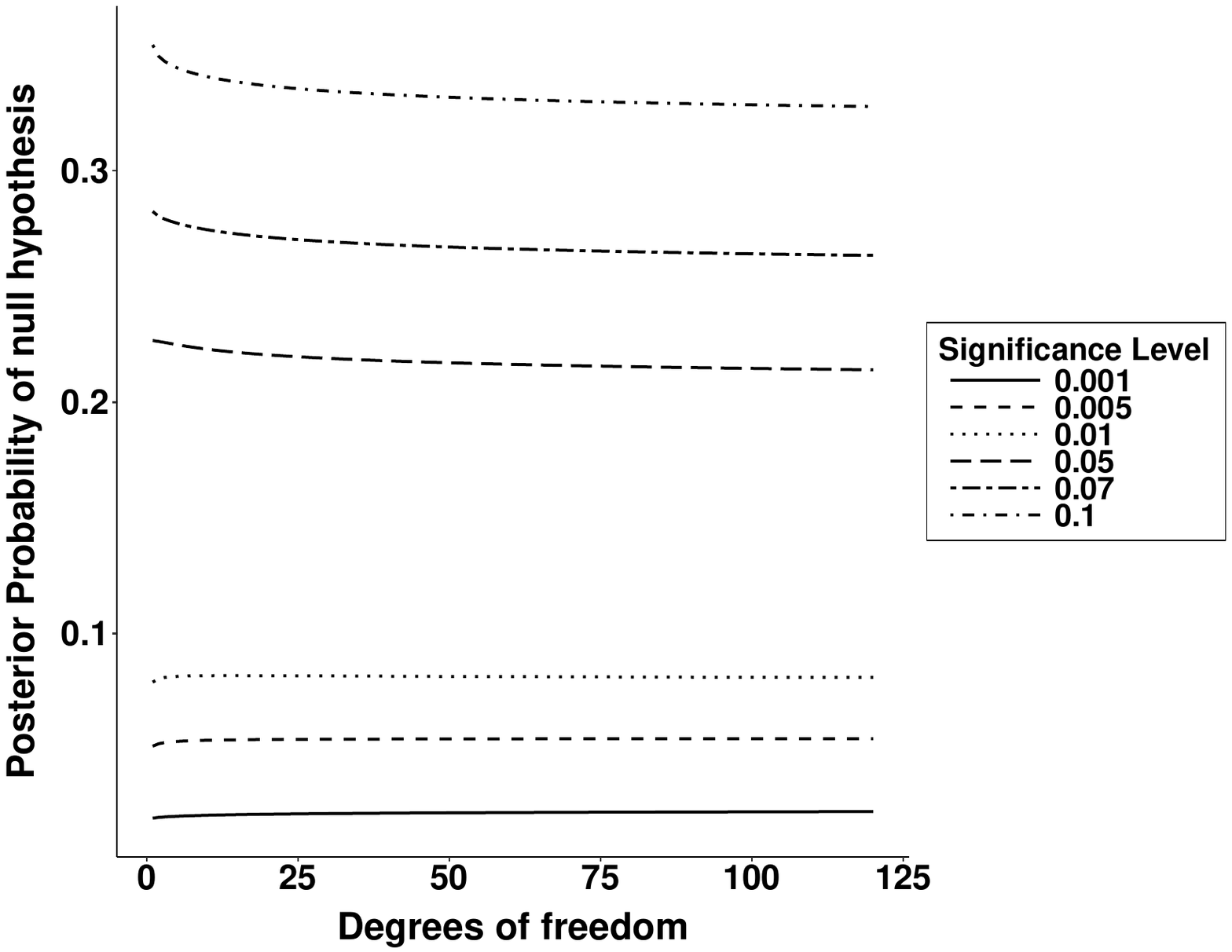}
		\caption{Posterior probability of the null hypothesis vs. degrees of freedom in chi-squared tests for different frequentist significance levels and $\gamma=3$. Based on the assignment of prior probability $1/2$ to both hypotheses. }
		\label{c5:ppnull}
	\end{figure}
	
	The posterior probabilities of null hypotheses with respect to different degrees of freedom in $\chi^2$ tests are depicted in Figure \ref{c5:ppnull}. These probabilities were computed under the assumption that the null and alternative hypotheses were equally likely {\em a priori}. The evidence threshold used to construct this plot was $\gamma=3.67$. 
	
	
	Figure \ref{c5:ppnull} shows that for $p$-values of $0.05$, the posterior probabilities assigned to the null were between 21.4\% and 22.7\%. The posterior probabilities resulting from $p$-values equal to $0.005$ were between 5.1\% and 5.5\%.
	
\begin{figure}[!t]
		\centering
		\includegraphics[width=75mm]{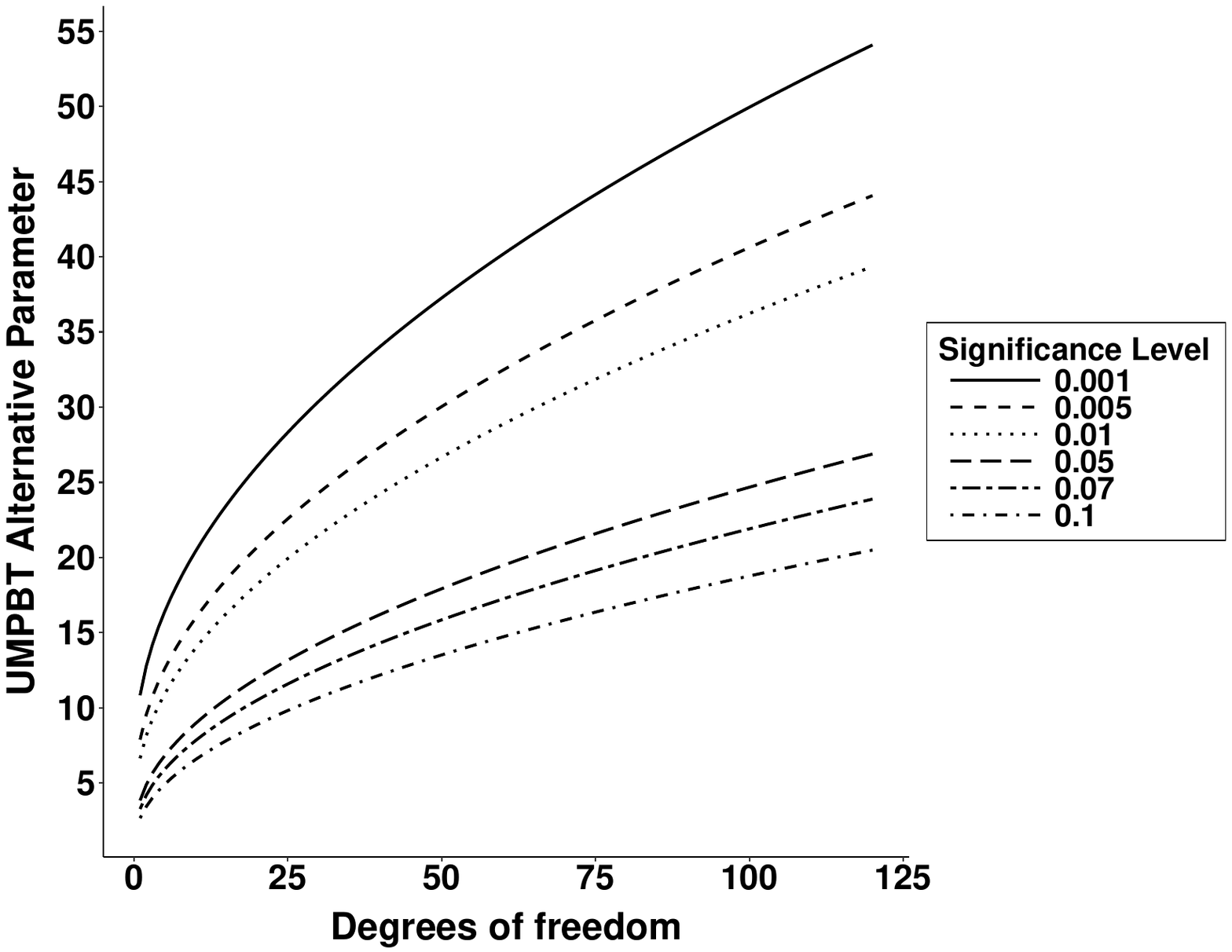}
		\caption{Alternative hypothesis values for UMPBT's matched to classical significance tests. }
		\label{fig3}
	\end{figure}

Figure \ref{fig3} depicts the UMPBT alternative hypothesis values in tests for which $\omgt$ has been matched to the rejection regions of classical chi-squared tests for various degrees of freedom and Type I errors.  In general, the value of the non-centrality parameter that defines the alternative hypothesis for a UMPBT must be determined numerically.  An R function that performs this calculation is provided in Supplementary Material.
	
	\subsection{One Sample t-test}\label{c5:onest}
	In this section we consider the one sample t-test for a normal mean when the variance, $\sigma^2$, is not known. We demonstrate that uniformly most powerful Bayesian tests do not exist in this setting. 
	
	Let $\bfx= \{x_1,x_2,\cdots,x_n\}$ represent $n$ i.i.d Gaussian observations, and define $\bar{x}$ to be the sample mean. The sample variance is defined in the usual way as $s^2=\sum_{i=1}^n(x_i-\bar{x})^2/(n-1)$.

	Suppose the prior distribution for $\sigma^2$ is inverse gamma with parameters $\alpha$ and $\beta$. Considering simple hypotheses, for every $\theta\in\Theta$ the marginal distribution of $\bfx$ is obtained by integrating out $\sigma^2$, leading to
	
	\begin{equation}\label{tmargin}
	m(\bfx) = (2\pi)^{n/2}\frac{\beta^\alpha\times\Gamma(n/2+\alpha)}{\Gamma(\alpha)\big[U+n(\bar{x}-\theta)^2\big]^{n/2+\alpha}}.
	\end{equation}
	Here, $U=\sum_{i=1}^n(x_i-\bar{x})^2+2\beta$. It follows that the Bayes factor for the test of $H_0: \theta = \theta_0$ versus $H_1:\theta = \theta_1$ can be expressed as
	
	\begin{equation}\label{tbf}
	g(y,\theta_1)= \Big[\frac{U+n(y-\theta_0)^2}{U+n(y-\theta_1)^2}\Big]^{n/2+\alpha},
	\end{equation}
	where $y=\bar{x}$.
	
	To use theorems exposed in the previous section, we first determine the form of $\Omega_\gamma(\theta)$. Letting $\gamma_n=\gamma^{\frac{2}{n+2\alpha}}$, it can be shown \citep{pnas} that for $\theta_1\in\Theta$, $\Omega_\gamma(\theta_1)$ can be expressed as
	
	\begin{equation}\label{trr}
	\Omega_\gamma(\theta_1) = \Big\{y:y^2+ \frac{2(\gamma_n\theta_1-\theta_0)}{1-\gamma_n}y-\big(\frac{\gamma_n\theta_1^2-\theta_0^2}{1-\gamma_n}-\frac{U}{n}\big)<0\Big\}.
	\end{equation}
	The roots of the quadratic function in (\ref{trr}) can be written as
	
	\begin{equation}\label{roots}
	a(\theta_1)=\frac{\gamma_n\theta_1-\theta_0}{\gamma_n-1}-\sqrt{\frac{\gamma_n(\theta_1-\theta_0)^2}{(\gamma_n-1)^2}-\frac{U}{n}}\quad \mbox{ and } \quad b(\theta_1)=\frac{\gamma_n\theta_1-\theta_0}{\gamma_n-1}+\sqrt{\frac{\gamma_n(\theta_1-\theta_0)^2}{(\gamma_n-1)^2}-\frac{U}{n}}.
	\end{equation}
	
	Thus,
	\begin{equation}\label{trr2}
	\Omega_\gamma(\theta_1)=\big\{y: a(\theta_1) < y < b(\theta_1)\big\}.
	\end{equation}
	
		From (\ref{roots}), it follows that for $\theta>\theta_0$,
		\begin{equation}\label{endbh1}
		\min\limits_\theta a(\theta) = \theta_0 + \sqrt{\frac{U(\gamma^*-1)}{n}}, \quad \argmin\limits_\theta a(\theta) =  \theta_0 + \sqrt{\frac{U(\gamma^*-1)}{n}},
		\end{equation}
		\begin{equation}
		\mbox{and } b(\theta) \mbox{ is a monotone increasing function of } \theta \mbox{ with } \lim\limits_{\theta\to +\infty} b(\theta) = +\infty.
		\end{equation}
		
		Similarly, for $\theta<\theta_0$,
		\begin{equation} 
		\max\limits_\theta b(\theta) = \theta_0 - \sqrt{\frac{U(\gamma^*-1)}{n}}, \quad \argmax\limits_\theta b(\theta) =  \theta_0 - \sqrt{\frac{U(\gamma^*-1)}{n}}, 
		\end{equation}
		\begin{equation}\label{endbh2}
		\mbox{and } a(\theta) \mbox{ is a monotone increasing function of } \theta \mbox{ with } \lim\limits_{\theta\to -\infty} a(\theta) = -\infty.
		\end{equation}
	
 It follows from (\ref{endbh1})--(\ref{endbh2}) that for a fixed $n < \infty$, no value of $\theta_1$ can be found to achieve the infimum of $a(\theta)$ and the supremum $b(\theta)$ at the same time, so that (\ref{c5:covereq}) cannot be achieved for any $\theta^*$. Indeed, different values of $\theta$ lead to non-nested $\omgt$, so that Corollary \ref{c5:sufthm} does not apply.
	
	To use Fact \ref{noexlem} to show that a UMPBT does not exist, consider two data generating parameters, say $\theta_{t}=2$ and $\theta_{t}=4$. Suppose the data-generating variance, $\sigma^2$, is equal to 1 and the evidence threshold, $\gamma$, is equal to 3. Take $\alpha=\beta=0$ so that a non-informative prior is assumed for $\sigma^2$. It follows from numerical analysis that the most powerful alternative when $\theta_t=2$ is obtained by taking $\theta^*=1.496$, while for $\theta_t=4$ the most powerful alternative is $\theta^*=2.394$.  Thus, a UMPBT does not exist for this test.
	
	Finally, we note that in classical hypothesis testing, \citet{lehmann} showed that $t$-tests for significance levels less than $0.5$ are not Uniformly Most Powerful in the classical sense. However, as noted by \citet{diaconis}, both one-sided and two-sided $t$-tests are UMP when attention is restricted to the class of all unbiased tests.
	
	\section{Discussion}\label{c5:discuss}
UMPBT's provide a new class of objective Bayesian hypothesis tests.  These tests facilitate a comparison between $p$-values from classical tests and Bayes factors from Bayesian tests.  These comparisons can be made by matching the rejection regions of the classical tests to the regions for which the Bayes factors from the UMPBT's exceed a specified evidence threshold.
	
  When UMPBT's are defined with evidence thresholds that remain fixed as sample sizes increase, they inherit certain inconsistencies of classical tests.  For instance, in large sample settings there remains a nonzero probability that the alternative hypothesis will be favored
		by the UMPBT even when the null hypothesis is true.  This deficiency stems from using a fixed evidence threshold (corresponding to a fixed significance level) as $n$ increases, which allows the alternative hypothesis defining the UMPBT to become arbitrarily close to the null hypothesis.  In general, we recommend increasing the evidence threshold with $n$ to avoid this problem.  
		Further discussion of this point is provided in \citet{umpbt}.   We also note that UMPBT's cannot generally be applied in sequential testing setting, since in most cases the alternative hypothesis that defines the UMPBT is determined by the sample size upon which the test is based.
	
	The primary contribution of this article is the extension of UMPBT's to a larger class of models and the introduction of a sufficient condition for the existence of UMPBT's. A practical mechanism for establishing the existence of a UMPBT was also provided. In cases when the sufficient condition is not satisfied, a procedure to verify that a UMPBT does not exist was provided by Fact \ref{noexlem}. These results allowed us to establish the existence of UMPBT's for tests of non-centrality parameters in $\chi^2$ statistics, which were illustrated for tests of independence in contingency tables.  By basing Bayes factors based on test statistics \citep{bfts}, the $\chi^2$ test can also be extended to obtain Bayes factors from likelihood ratio tests and score tests, which are among the most commonly used classical test statistics.  We also showed that uniformly most powerful Bayesian tests do not exist for one-sample $t$-tests.
	
 It is important to note that UMPBT's do not provide an upper bound on the Bayes factor against a point null hypothesis, a common misperception.  When the null hypothesis is true (i.e., represents the data-generating parameter), the Bayes factor against the null hypothesis will typically be substantially smaller than the likelihood ratio statistic for the test, which provides an actual upper bound on the Bayes factor.  Conversely when the alternative hypothesis is true, the Bayes factor in favor of the alternative hypothesis will typically be smaller than the likelihood ratio statistic, particularly if the data-generating parameter exceeds the value defined under the UMPBT alternative hypothesis.  It is only when the maximum likelihood estimate of the tested parameter is close to a (simple) UMPBT alternative that the Bayes factor based on the UMPBT provides an approximate upper bound on the Bayes factor.
		
		To illustrate the importance of this difference, suppose $\bar{x} = 0$ in a test of a normal mean in which the UMPBT's alternative hypothesis is chosen so that $\omgts$ matches the one-sided 5\% classical test's rejection region.  Then the Bayes factor based on the UMPBT is $BF_{10}(\bar{x}) = .258$, whereas the maximum Bayes factor is 1.  In this case, the use of the UMPBT results in a posterior probability for the null hypothesis that can significantly exceed 0.5 (i.e., 0.795) when both hypotheses are assigned equal prior probability, something that is not possible if the MLE is used to set the alternative hypothesis. Similarly, if $\bar{x}=2*1.645 \sigma/\sqrt{n}$, or twice the UMPBT alternative, then $BF_{10} = 57.9$, whereas the maximum Bayes factor is 224.1.  That is, the UMPBT-based Bayes factor is a factor of 4 smaller than the likelihood ratio statistic.

	In future research, we hope to extend UMPBT's to Bayesian variable selection problems and examine constraints that will allow this methodology to be extended to multi-dimensional exponential family distributions. 

\bibliography{ref}
\end{document}